\documentclass[11pt]{article}
\usepackage[lmargin=1in,rmargin=1in]{geometry}
\usepackage{amsmath}
\usepackage{amssymb}
\usepackage{listings}

\newcommand{\ssage}{{\em Sage\ }}

\usepackage{url}

\title{Program for calculating bounds on the minimum rank of a graph using Sage}

\author{Laura DeLoss\thanks{Department of Mathematics, Iowa State
    University, Ames, IA 50011, USA (delolau@iastate.edu,
    grout@iastate.edu, tmckay16@iastate.edu, smithj@iastate.edu,
    gtims@iastate.edu).} \and Jason Grout\footnotemark[1] \and Tracy
  McKay\footnotemark[1] \and Jason Smith\footnotemark[1]\and Geoff
  Tims\footnotemark[1]}

\begin{document}
\maketitle

\begin{abstract}
  The minimum rank of a simple graph $G$ is defined to be the smallest
  possible rank over all symmetric real matrices whose $ij$th entry
  (for $i\neq j$) is nonzero whenever $\{i,j\}$ is an edge in $G$ and is
  zero otherwise.  Minimum rank is a difficult parameter to compute.
  However, there are now a number of known reduction techniques and
  bounds that can be programmed on a computer; we have developed a
  program using the open-source mathematics software \ssage to
  implement several techniques.  In this note, we provide the source code
  for this program.
\end{abstract}


\bigskip

\section{Introduction}

In this note, we provide a listing for our \ssage
\cite[version~3.1.2]{sage} program which computes upper and lower
bounds for the minimum rank of a graph.

We first include several example \ssage sessions, illustrating the
main function of the program, \verb|minrank_bounds|.  The lines
starting with \verb|sage:| are input lines (the ``sage:'' is not
typed), while the other lines are output lines.  Without any options,
the \verb|minrank_bounds| function returns two numbers: a lower bound
and an upper bound for the minimum rank.  In the following example, we
find that $K_3$ has minimum rank 1 (i.e., the lower bound and the
upper bound are both 1).

\begin{verbatim}
sage: minrank_bounds(graphs.CompleteGraph(3)) 
(1, 1)
\end{verbatim}

With the \verb|all_bounds| option set to \verb|True|, the function
returns two listings in the format ``name: value'', where the name is
the name of a bound and the value is the upper or lower bound.  In the
following example, we calculate bounds on the minimum rank of the
Petersen graph.  In the second call, we see that lower bounds are
minimally given by the obvious bound on ranks (i.e., the rank has to
be at least zero) and also by a minimal zero forcing set (a lower
bound of 5).  For upper bounds, we have a clique cover which gives a
bound of 15, the fact that the graph is not outer planar (upper bound
of 7), the fact that the graph is not a path (upper bound of 8), the
fact that the graph is not planar (upper bound of 6), the fact that
the minimum rank is at most $|G|-1$ (upper bound of 9), and the
trivial upper bound on ranks (i.e., the rank is at most $|G|$).

\begin{verbatim}
sage: minrank_bounds(graphs.PetersenGraph())
(5, 6)
sage: minrank_bounds(graphs.PetersenGraph(), all_bounds=True)
({'rank': 0, 'zero forcing': 5},
{'clique cover': 15,
'not outer planar': 7,
'not path': 8,
'not planar': 6,
'order': 9,
'rank': 10})
\end{verbatim}

If the \verb|tests| option is set to a list, then only those tests are
run.  In the following example, we compute bounds on the minimum rank
of the Heawood graph by using the zero forcing test and by testing if
the graph is not planar.  The ``rank'' tests are always run (and
always give 0 for a lower bound and $|G|$ for an upper bound).

\begin{verbatim}
sage: minrank_bounds(graphs.HeawoodGraph(), all_bounds=True, 
sage:   tests=['zero forcing', 'not planar'])
({'rank': 0, 'zero forcing': 8}, {'not planar': 10, 'rank': 14})
\end{verbatim}

For more complete documentation and a list of tests that can be run,
as well as for several more examples, print the help by typing the
function name followed by a question mark: \verb|minrank_bounds?|.

This program contains and uses the minimum rank data listed in
\cite{SMALLGRAPHS}.

\section{Program listing}

In the typeset program listing below, some lines are automatically
broken that are not actually broken in the source code.  If a line is
broken into two lines in the listing, but should appear as one line in
the program, then the line will end with $\hookleftarrow$ and the remainder of the
line will start with $\hookrightarrow$.

To use this program, download and extract the ``source'' from
arxiv.org.  Then you can
\begin{itemize}
\item upload the accompanying \verb|minrank.sws| \ssage worksheet into
  a notebook, or
\item copy the contents of the accompanying \verb|minrank.sage| file
  into a cell of a \ssage worksheet, or
\item load the \verb|minrank.sage| file into a running \ssage terminal
  session.
\end{itemize}

We now proceed with the program listing.

\begin{lstlisting}
##########################################################
# Imports all the graphs of order 7 or less and stores   #
# them in a list called atlas_graphs so that             #
# atlas_graphs[i] is the ith graph in the atlas of       #
# graphs                                                 #
##########################################################


import networkx.generators.atlas
atlas_graphs = [Graph(i) for i in \
                    networkx.generators.atlas.graph_atlas_g()]

##########################################################
# A list "database" of all the minimum ranks of graphs   #
# of order 7 or less.                                    #
#                                                        #
# The minimum ranks are stored as ordered pairs: the     #
# first coordinate is the graph number, the second       #
# coordinate is the minimum rank.  The first tuple in    #
# the list min_ranks is just a position holder.          #
##########################################################

min_ranks = [(0,None), (1,0), (2,0), (3,1), (4,0), (5,1), (6,2),
(7,1), (8,0), (9,1), (10,2), (11,2), (12,1), (13,2), (14,3), (15,2),
(16,2), (17,2), (18,1), (19,0), (20,1), (21,2), (22,2), (23,1),
(24,2), (25,3), (26,3), (27,2), (28,2), (29,2), (30,3), (31,4),
(32,2), (33,2), (34,3), (35,3), (36,3), (37,3), (38,3), (39,1),
(40,3), (41,3), (42,2), (43,3), (44,2), (45,2), (46,2), (47,3),
(48,2), (49,2), (50,2), (51,2), (52,1), (53,0), (54,1), (55,2),
(56,2), (57,1), (58,2), (59,3), (60,3), (61,3), (62,2), (63,2),
(64,2), (65,3), (66,4), (67,2), (68,3), (69,4), (70,4), (71,2),
(72,3), (73,3), (74,3), (75,3), (76,3), (77,2), (78,3), (79,4),
(80,4), (81,4), (82,3), (83,5), (84,3), (85,3), (86,1), (87,3),
(88,3), (89,2), (90,3), (91,2), (92,3), (93,4), (94,4), (95,4),
(96,4), (97,4), (98,4), (99,4), (100,3), (101,3), (102,4), (103,4),
(104,4), (105,4), (106,2), (107,2), (108,2), (109,3), (110,2),
(111,4), (112,4), (113,4), (114,3), (115,4), (116,2), (117,3),
(118,4), (119,3), (120,4), (121,3), (122,4), (123,4), (124,4),
(125,3), (126,3), (127,4), (128,4), (129,3), (130,3), (131,2),
(132,2), (133,3), (134,3), (135,3), (136,4), (137,4), (138,3),
(139,4), (140,3), (141,3), (142,3), (143,3), (144,3), (145,3),
(146,2), (147,4), (148,4), (149,3), (150,3), (151,3), (152,4),
(153,3), (154,3), (155,2), (156,3), (157,3), (158,3), (159,3),
(160,3), (161,2), (162,3), (163,3), (164,4), (165,2), (166,3),
(167,4), (168,3), (169,3), (170,3), (171,3), (172,3), (173,3),
(174,3), (175,2), (176,1), (177,3), (178,3), (179,3), (180,3),
(181,3), (182,3), (183,3), (184,3), (185,3), (186,3), (187,3),
(188,3), (189,2), (190,2), (191,2), (192,3), (193,3), (194,2),
(195,2), (196,3), (197,2), (198,3), (199,2), (200,2), (201,2),
(202,3), (203,2), (204,2), (205,2), (206,2), (207,2), (208,1),
(209,0), (210,1), (211,2), (212,2), (213,1), (214,2), (215,3),
(216,3), (217,3), (218,2), (219,2), (220,2), (221,3), (222,4),
(223,2), (224,3), (225,4), (226,4), (227,4), (228,2), (229,3),
(230,3), (231,3), (232,3), (233,3), (234,2), (235,3), (236,4),
(237,4), (238,4), (239,3), (240,5), (241,3), (242,3), (243,3),
(244,4), (245,4), (246,5), (247,5), (248,3), (249,1), (250,3),
(251,3), (252,2), (253,3), (254,2), (255,3), (256,4), (257,4),
(258,4), (259,4), (260,4), (261,4), (262,3), (263,4), (264,3),
(265,4), (266,4), (267,4), (268,4), (269,2), (270,2), (271,3),
(272,4), (273,4), (274,4), (275,4), (276,5), (277,4), (278,4),
(279,5), (280,5), (281,4), (282,4), (283,4), (284,5), (285,3),
(286,6), (287,4), (288,4), (289,4), (290,2), (291,2), (292,3),
(293,2), (294,4), (295,4), (296,4), (297,3), (298,4), (299,2),
(300,3), (301,4), (302,3), (303,4), (304,3), (305,4), (306,4),
(307,4), (308,3), (309,3), (310,4), (311,3), (312,4), (313,3),
(314,3), (315,4), (316,4), (317,5), (318,4), (319,4), (320,5),
(321,5), (322,5), (323,4), (324,4), (325,4), (326,3), (327,5),
(328,5), (329,4), (330,4), (331,5), (332,5), (333,5), (334,5),
(335,3), (336,5), (337,5), (338,5), (339,4), (340,5), (341,5),
(342,5), (343,4), (344,4), (345,4), (346,4), (347,3), (348,5),
(349,5), (350,5), (351,5), (352,3), (353,5), (354,3), (355,2),
(356,2), (357,3), (358,3), (359,3), (360,4), (361,4), (362,3),
(363,4), (364,3), (365,3), (366,3), (367,3), (368,3), (369,3),
(370,2), (371,4), (372,4), (373,3), (374,3), (375,3), (376,3),
(377,4), (378,3), (379,4), (380,4), (381,5), (382,3), (383,5),
(384,4), (385,5), (386,4), (387,3), (388,4), (389,4), (390,5),
(391,5), (392,4), (393,5), (394,5), (395,4), (396,4), (397,3),
(398,5), (399,5), (400,5), (401,5), (402,5), (403,4), (404,4),
(405,4), (406,4), (407,4), (408,4), (409,4), (410,4), (411,4),
(412,5), (413,5), (414,5), (415,4), (416,4), (417,3), (418,3),
(419,4), (420,4), (421,5), (422,5), (423,5), (424,4), (425,4),
(426,4), (427,5), (428,4), (429,4), (430,4), (431,4), (432,5),
(433,5), (434,5), (435,5), (436,4), (437,5), (438,5), (439,5),
(440,4), (441,4), (442,4), (443,4), (444,4), (445,5), (446,5),
(447,4), (448,4), (449,4), (450,4), (451,3), (452,2), (453,3),
(454,3), (455,3), (456,3), (457,3), (458,2), (459,3), (460,3),
(461,4), (462,2), (463,3), (464,4), (465,3), (466,3), (467,3),
(468,3), (469,3), (470,3), (471,3), (472,2), (473,3), (474,4),
(475,4), (476,4), (477,4), (478,5), (479,5), (480,4), (481,4),
(482,5), (483,4), (484,5), (485,4), (486,4), (487,4), (488,5),
(489,5), (490,4), (491,4), (492,4), (493,4), (494,4), (495,4),
(496,3), (497,5), (498,4), (499,4), (500,4), (501,3), (502,3),
(503,4), (504,4), (505,4), (506,4), (507,3), (508,5), (509,5),
(510,4), (511,4), (512,5), (513,4), (514,4), (515,5), (516,5),
(517,5), (518,5), (519,4), (520,4), (521,4), (522,4), (523,4),
(524,4), (525,3), (526,5), (527,5), (528,5), (529,5), (530,5),
(531,4), (532,4), (533,5), (534,4), (535,4), (536,4), (537,4),
(538,4), (539,4), (540,4), (541,4), (542,4), (543,4), (544,4),
(545,4), (546,4), (547,4), (548,5), (549,4), (550,4), (551,3),
(552,4), (553,4), (554,3), (555,4), (556,4), (557,3), (558,3),
(559,5), (560,4), (561,5), (562,4), (563,5), (564,4), (565,4),
(566,5), (567,4), (568,4), (569,4), (570,3), (571,4), (572,4),
(573,4), (574,5), (575,5), (576,4), (577,4), (578,4), (579,4),
(580,4), (581,4), (582,2), (583,1), (584,3), (585,3), (586,3),
(587,3), (588,3), (589,3), (590,3), (591,3), (592,3), (593,3),
(594,3), (595,3), (596,2), (597,2), (598,4), (599,4), (600,4),
(601,4), (602,4), (603,4), (604,4), (605,4), (606,4), (607,4),
(608,4), (609,4), (610,3), (611,3), (612,3), (613,4), (614,4),
(615,3), (616,4), (617,4), (618,5), (619,3), (620,4), (621,4),
(622,5), (623,5), (624,4), (625,4), (626,4), (627,4), (628,4),
(629,4), (630,3), (631,4), (632,5), (633,4), (634,4), (635,4),
(636,4), (637,4), (638,4), (639,4), (640,5), (641,4), (642,4),
(643,4), (644,3), (645,4), (646,5), (647,4), (648,4), (649,4),
(650,4), (651,4), (652,4), (653,4), (654,4), (655,4), (656,4),
(657,4), (658,4), (659,4), (660,4), (661,4), (662,4), (663,4),
(664,4), (665,4), (666,4), (667,3), (668,3), (669,3), (670,2),
(671,4), (672,4), (673,4), (674,4), (675,4), (676,4), (677,4),
(678,3), (679,4), (680,4), (681,3), (682,5), (683,5), (684,4),
(685,4), (686,3), (687,3), (688,4), (689,3), (690,5), (691,4),
(692,4), (693,4), (694,5), (695,4), (696,4), (697,4), (698,4),
(699,4), (700,4), (701,4), (702,4), (703,4), (704,4), (705,4),
(706,4), (707,4), (708,4), (709,4), (710,5), (711,4), (712,4),
(713,4), (714,4), (715,4), (716,4), (717,4), (718,4), (719,4),
(720,4), (721,3), (722,4), (723,4), (724,4), (725,3), (726,3),
(727,4), (728,4), (729,4), (730,3), (731,2), (732,3), (733,3),
(734,2), (735,2), (736,3), (737,2), (738,3), (739,2), (740,4),
(741,4), (742,4), (743,3), (744,4), (745,2), (746,4), (747,4),
(748,4), (749,4), (750,4), (751,4), (752,4), (753,4), (754,4),
(755,4), (756,4), (757,4), (758,4), (759,4), (760,4), (761,4),
(762,4), (763,4), (764,4), (765,4), (766,4), (767,4), (768,4),
(769,4), (770,4), (771,4), (772,4), (773,4), (774,4), (775,4),
(776,4), (777,4), (778,4), (779,4), (780,3), (781,3), (782,4),
(783,4), (784,4), (785,4), (786,3), (787,3), (788,4), (789,3),
(790,2), (791,3), (792,4), (793,4), (794,4), (795,4), (796,3),
(797,4), (798,4), (799,4), (800,4), (801,3), (802,4), (803,4),
(804,4), (805,4), (806,4), (807,4), (808,4), (809,3), (810,4),
(811,4), (812,3), (813,5), (814,3), (815,3), (816,4), (817,4),
(818,4), (819,4), (820,4), (821,5), (822,4), (823,4), (824,4),
(825,4), (826,4), (827,4), (828,4), (829,3), (830,4), (831,3),
(832,3), (833,4), (834,4), (835,4), (836,4), (837,4), (838,4),
(839,4), (840,5), (841,4), (842,4), (843,4), (844,4), (845,4),
(846,3), (847,4), (848,4), (849,4), (850,4), (851,3), (852,4),
(853,4), (854,4), (855,4), (856,3), (857,4), (858,4), (859,4),
(860,4), (861,4), (862,4), (863,3), (864,4), (865,3), (866,4),
(867,4), (868,4), (869,3), (870,4), (871,4), (872,3), (873,3),
(874,3), (875,4), (876,3), (877,4), (878,3), (879,2), (880,2),
(881,3), (882,2), (883,2), (884,3), (885,3), (886,4), (887,4),
(888,4), (889,4), (890,4), (891,4), (892,4), (893,3), (894,3),
(895,3), (896,3), (897,4), (898,3), (899,4), (900,4), (901,3),
(902,3), (903,4), (904,4), (905,4), (906,3), (907,4), (908,4),
(909,3), (910,4), (911,3), (912,3), (913,3), (914,4), (915,4),
(916,4), (917,4), (918,3), (919,4), (920,4), (921,4), (922,4),
(923,4), (924,3), (925,3), (926,4), (927,4), (928,4), (929,4),
(930,4), (931,3), (932,3), (933,4), (934,4), (935,4), (936,4),
(937,4), (938,4), (939,4), (940,4), (941,4), (942,4), (943,4),
(944,3), (945,3), (946,4), (947,3), (948,3), (949,3), (950,4),
(951,4), (952,4), (953,3), (954,4), (955,4), (956,3), (957,3),
(958,3), (959,4), (960,4), (961,4), (962,4), (963,4), (964,4),
(965,4), (966,4), (967,4), (968,4), (969,4), (970,3), (971,4),
(972,4), (973,3), (974,4), (975,3), (976,4), (977,3), (978,3),
(979,4), (980,4), (981,4), (982,4), (983,3), (984,3), (985,4),
(986,4), (987,3), (988,3), (989,4), (990,3), (991,3), (992,4),
(993,4), (994,3), (995,3), (996,3), (997,4), (998,4), (999,4),
(1000,3), (1001,3), (1002,3), (1003,3), (1004,3), (1005,3), (1006,4),
(1007,2), (1008,4), (1009,2), (1010,2), (1011,2), (1012,3), (1013,3),
(1014,3), (1015,4), (1016,4), (1017,3), (1018,3), (1019,3), (1020,3),
(1021,4), (1022,3), (1023,3), (1024,3), (1025,4), (1026,4), (1027,4),
(1028,3), (1029,4), (1030,4), (1031,4), (1032,2), (1033,4), (1034,3),
(1035,3), (1036,3), (1037,3), (1038,3), (1039,4), (1040,3), (1041,4),
(1042,3), (1043,4), (1044,3), (1045,3), (1046,4), (1047,4), (1048,4),
(1049,3), (1050,4), (1051,4), (1052,4), (1053,4), (1054,4), (1055,4),
(1056,3), (1057,3), (1058,4), (1059,4), (1060,3), (1061,4), (1062,3),
(1063,3), (1064,3), (1065,4), (1066,3), (1067,3), (1068,3), (1069,4),
(1070,3), (1071,4), (1072,3), (1073,3), (1074,3), (1075,3), (1076,3),
(1077,3), (1078,4), (1079,3), (1080,4), (1081,3), (1082,4), (1083,4),
(1084,3), (1085,3), (1086,3), (1087,3), (1088,2), (1089,4), (1090,3),
(1091,4), (1092,3), (1093,4), (1094,3), (1095,3), (1096,3), (1097,4),
(1098,3), (1099,3), (1100,3), (1101,4), (1102,3), (1103,3), (1104,3),
(1105,3), (1106,3), (1107,2), (1108,3), (1109,3), (1110,3), (1111,3),
(1112,3), (1113,3), (1114,3), (1115,3), (1116,3), (1117,4), (1118,4),
(1119,3), (1120,3), (1121,4), (1122,3), (1123,3), (1124,3), (1125,3),
(1126,3), (1127,4), (1128,3), (1129,3), (1130,3), (1131,3), (1132,3),
(1133,3), (1134,3), (1135,3), (1136,3), (1137,3), (1138,3), (1139,3),
(1140,2), (1141,4), (1142,4), (1143,3), (1144,3), (1145,4), (1146,3),
(1147,3), (1148,3), (1149,3), (1150,4), (1151,3), (1152,3), (1153,3),
(1154,4), (1155,3), (1156,3), (1157,3), (1158,3), (1159,3), (1160,4),
(1161,3), (1162,3), (1163,3), (1164,2), (1165,3), (1166,3), (1167,3),
(1168,3), (1169,3), (1170,3), (1171,2), (1172,1), (1173,3), (1174,3),
(1175,3), (1176,3), (1177,3), (1178,3), (1179,3), (1180,3), (1181,3),
(1182,3), (1183,3), (1184,2), (1185,3), (1186,3), (1187,4), (1188,2),
(1189,3), (1190,4), (1191,3), (1192,3), (1193,3), (1194,3), (1195,3),
(1196,3), (1197,3), (1198,3), (1199,3), (1200,3), (1201,3), (1202,3),
(1203,3), (1204,3), (1205,3), (1206,2), (1207,3), (1208,2), (1209,3),
(1210,3), (1211,2), (1212,3), (1213,2), (1214,3), (1215,3), (1216,2),
(1217,3), (1218,3), (1219,3), (1220,3), (1221,3), (1222,3), (1223,3),
(1224,3), (1225,3), (1226,2), (1227,3), (1228,3), (1229,2), (1230,2),
(1231,3), (1232,3), (1233,2), (1234,2), (1235,3), (1236,3), (1237,2),
(1238,2), (1239,3), (1240,2), (1241,3), (1242,2), (1243,2), (1244,2),
(1245,2), (1246,3), (1247,2), (1248,2), (1249,2), (1250,2), (1251,2),
(1252,1)]

###########################################################

# The atlas of graphs is ordered by number of vertices, then number of
# edges.  To make our search for minimum ranks more efficient, we
# store the first and last indices of the graphs having a certain
# number of vertices and edges.

# graph_help[n][m] is a tuple (i,j), meaning that the graphs with n
# vertices and m edges start at index i and end on index j in the
# atlas of graphs.
###########################################################



graph_help = {
 1: {0: (1,1)},

 2: {0: (2,2), 1: (3,3)},

 3: {0: (4,4), 1: (5,5), 2: (6,6), 3: (7,7)},

 4: {0: (8,8), 1: (9,9), 2: (10,11), 3: (12,14), 4: (15,16), 
     5: (17,17), 6: (18,18)},
 
 5: {0: (19,19), 1: (20,20), 2: (21,22), 3: (23,26), 4: (27,32), 5:
         (33,38), 6: (39,44), 7: (45,48), 8: (49,50), 9: (51,51), 10:
         (52,52)},
 
 6: {0: (53,53), 1:(54,54), 2:(55,56), 3: (57,61), 4: (62,70), 
     5: (71,85), 6: (86,106), 7: (107,130), 8: (131,154), 9: (155,175), 
     10: (176,190), 11: (191,199), 12: (200,204), 13: (205,206), 
     14: (207,207), 15: (208,208)},

 7: {0: (209,209), 1: (210,210), 2: (211,212), 3: (213,217), 
     4: (218,227), 5: (228,248), 6: (249,289), 7: (290,354), 8: (355,451), 
     9: (452,582), 10: (583,730), 11: (731,878), 12: (879,1009), 
     13: (1010,1106), 14: (1107,1171), 15: (1172,1212), 16: (1213,1233), 
     17: (1234,1243), 18: (1244,1248), 19: (1249,1250), 20: (1251,1251), 
     21: (1252,1252)}}


def get_mr_from_list(graph):
    """
    Return the minimum rank of a graph of order 7 or less from the
    list min_ranks.

    INPUT:
        graph -- the graph whose minimum rank is to be 
        found

    OUTPUT:
        the minimum rank if the graph is in the list 
        or False if it is not

    EXAMPLES:
        sage: get_mr_from_list(Graph({0:[2,3],1:[2],2:[3,4],3:[4]}))
        3
        sage: get_mr_from_list(graphs.PathGraph(8))
        False
    
    """

    #check to make sure graph can be found in list
    if graph.order()>7:
        return False
    
    order = graph.order() #number of vertices
    size = graph.size() #number of edges

    starting_index, ending_index = graph_help[order][size]
    import pdb
    #look for graph and return minimum rank
    for i in [starting_index..ending_index]:
        if graph.is_isomorphic(atlas_graphs[i]):
            return min_ranks[i][1]

    raise ValueError, "This should never happen!"


def zerosgame(graph, initial_set=[]):
   """
   Apply the color-change rule to a given graph given an optional
   initial set.

   INPUT:
        graph -- the graph on which to apply the rule
        initial_set -- the set of "zero" (black) vertices in the graph

   OUTPUT:
        the list of zero (black) vertices in the resulting derived
        coloring

   EXAMPLES:
        sage: zerosgame(graphs.PathGraph(5))
        []
        sage: zerosgame(graphs.PathGraph(5),[0])
        [0, 1, 2, 3, 4]
   """       

   new_zero_set=set(initial_set)
   zero_set=set([])
   zero_neighbors={}
   active_zero_set = set([])
   inactive_zero_set = set([])
   another_run=True
   while another_run:
       another_run=False
       # Add the new zero vertices
       zero_set.update(new_zero_set)
       active_zero_set.update(new_zero_set)
       active_zero_set.difference_update(inactive_zero_set)
       zero_neighbors.update([[i, 
                   set(graph.neighbors(i)).difference(zero_set)] 
                              for i in new_zero_set])
       # Find the next set of zero vertices
       new_zero_set.clear()
       inactive_zero_set.clear()
       for v in active_zero_set:
           zero_neighbors[v].difference_update(zero_set)
           if len(zero_neighbors[v])==1:
               new_zero_set.add(zero_neighbors[v].pop())
               inactive_zero_set.add(v)
               another_run=True
   return list(zero_set)


def find_zero_forcing_set(graph, bound=None):
   """
   Return a zero forcing set of minimum order that also has order
   less than the given bound.
    
   INPUT:
      graph -- the graph on which to find the zero-forcing set
      bound -- the maximum acceptable order for a zero-forcing set

   OUTPUT:
      a zero-forcing set of minimum order that also has order less
      than the bound if one exists; False if no such zero-forcing set
      can be found

   EXAMPLES:
      sage: find_zero_forcing_set(graphs.CompleteGraph(5))
      {0, 1, 2, 3}
      sage: find_zero_forcing_set(graphs.CompleteGraph(5),2)
      False
   """
   order=graph.order()
   if bound is None:
       bound = order
   if bound<0:
       bound=1
   vertices=graph.vertices()
   mindegree=min(graph.degree())
   for i in [mindegree..bound]:
       for subset in Subsets(vertices,i):
           outcome=zerosgame(graph,subset)
           if len(outcome)==order:
               return subset
   return False


def find_Z(graph):
    """
    Returns the order of a smallest zero-forcing set of a graph

    INPUT:
        graph -- the graph on which to find a smallest zero-forcing set

    OUTPUT:
        the minimum possible order for a zero-forcing set of the graph

    EXAMPLES:
        sage: find_Z(graphs.CompleteGraph(5))
        4      

    """
    return len(find_zero_forcing_set(graph))


def has_forbidden_induced_subgraph(graph):
    """
    Check for a forbidden induced subgraph (a path on 4 vertices,
    fish, dart, or complete tripartite graph).

    INPUT:
        graph -- the graph to be checked 

    OUTPUT:
        True if the graph contains an induced copy of P_4, fish, dart,
        or K_{3,3,3}; False if it does not

    EXAMPLES:
        sage: has_forbidden_induced_subgraph(graphs.CompleteGraph(10))
        False
        sage: has_forbidden_induced_subgraph(graphs.PathGraph(3))
        False
        sage: K333 = Graph({0: [3,4,5,6,7,8], 1: [3,4,5,6,7,8], 2: [3,4,5,6,7,8], 3: [6,7,8], 4: [6,7,8], 5: [6,7,8]})
        sage: has_forbidden_induced_subgraph(K333)
        True
        sage: g = Graph({0:[1,2,3,4], 1: [2]}) # fish
        sage: has_forbidden_induced_subgraph(g)
        True
        sage: g.add_edge((0,6))
        sage: has_forbidden_induced_subgraph(g)
        True
    """
    order = graph.order()
    vertices = graph.vertices()
    path = atlas_graphs[14]
    fish = atlas_graphs[34]
    dart = atlas_graphs[40]
    K333 = Graph({0: [3,4,5,6,7,8], 1: [3,4,5,6,7,8], \
                  2: [3,4,5,6,7,8], 3: [6,7,8], \
                  4: [6,7,8], 5: [6,7,8]})
    if order < 4:
        return False
    for sub_vertices in Combinations(vertices,4): 
        # Finds all order 4 induced subgraphs
        if graph.subgraph(sub_vertices).is_isomorphic(path):
            return True
    if order < 5:
        return False
    for sub_vertices in Combinations(vertices,5): 
        # Finds all order 5 induced subgraphs
        if graph.subgraph(sub_vertices).is_isomorphic(dart):
            return True
        if graph.subgraph(sub_vertices).is_isomorphic(fish):
            return True
    if order < 9:
        return False
    for sub_vertices in Combinations(vertices,9): 
        # Finds all order 9 induced subgraphs
        if graph.subgraph(sub_vertices).is_isomorphic(K333):
            return True
    return False


def is_outerplanar(graph):
    """
    Check if the graph is outer-planar

    INPUT:
        graph -- is the graph to be checked

    OUTPUT:
        True if the graph is outer-planar; False if it is not

    EXAMPLES:
        sage: is_outerplanar(Graph({0:[1,2,3,4],1:[2,4],2:[3],3:[4]}))
        False 
        sage: is_outerplanar(graphs.CompleteGraph(3))
        True            
    """
    h = graph.copy()
    h.delete_vertices([v for v in h.vertices() if h.degree(v) == 0])
    if h.order()==0:
        return True
    h.set_boundary(h.vertices())
    return h.is_circular_planar(ordered=False)



def find_edge_max_clique(edge,graph):
    """
    Given an edge and a graph, return a maximal clique of the graph
    that contains the edge.

    INPUT:
        edge -- the edge for which to find the maximal clique
        graph -- the graph containing the edge

    OUTPUT:
        the list of vertices in a maximal clique that contains the edge
        if the edge is not in the graph, it will return None

    EXAMPLES:
        sage: find_edge_max_clique((1,2), graphs.CompleteGraph(5))
        [0, 1, 2, 3, 4] 
        sage: find_edge_max_clique((1,3), graphs.PathGraph(5)) 
    """
    vertex1=edge[0] # first vertex of edge
    vertex2=edge[1] # second vertex of edge
    pot_cliques=graph.cliques_containing_vertex(vertex1)
    
    # sort the cliques containing vertex1 by order, largest first
    pot_cliques.sort(key=len, reverse=True)

    for clique in pot_cliques:
        if vertex2 in clique:
            return clique
    return None


def edge_clique_cover(graph, bound=None):
    """
    Assuming the graph is connected, return an edge clique cover for
    the graph if the number of covering cliques is at most bound;
    otherwise, returns None.

    INPUT:
        graph -- the graph
        bound -- the maximum number of cliques to consider

    OUTPUT:
        If a clique cover is found that has at most bound cliques, the
        clique cover is returned as a list of lists, each sublist
        being the vertices of a clique.
        
        If a clique cover from this function requires more than bound
        cliques, None is returned.

    EXAMPLES:
        sage: edge_clique_cover(graphs.PathGraph(3))
        [[0, 1], [1, 2]]
        sage: edge_clique_cover(graphs.CompleteGraph(5))
        [[0, 1, 2, 3, 4]]
        sage: edge_clique_cover(graphs.HouseGraph())
        [[2, 3, 4], [0, 1], [0, 2], [1, 3]]
        sage: edge_clique_cover(graphs.PetersenGraph(),bound=4)
        sage:
    """
    
    # Take care of trivial case
    if graph.size() == 0:
        return []

    max_cliques=graph.cliques()
    max_cliques.sort(key=len)
    largest_clique_vertices = len(max_cliques[-1])
    max_cliques = [sorted(clique) for clique in max_cliques]
    largest_clique_edges = largest_clique_vertices \
                           *(largest_clique_vertices-1)/2
    edges_of_graph=graph.edges(labels=False)
    num_edges = graph.size()

    mandatory_cliques=[]
    

    for v in graph.vertices():
        # If v is contained in only one clique, then that clique must
        # be in the clique cover
        cliques_containing_v = [c for c in max_cliques if v in c]
        if len(cliques_containing_v)==1 \
                and (cliques_containing_v[0] not in mandatory_cliques):
            mandatory_cliques.append(cliques_containing_v[0])
    for e in graph.edges():
        # If e is contained in only one clique, then that clique must
        # be in the clique cover
        cliques_containing_e = [c for c in max_cliques 
                                if e[0] in c and e[1] in c]
        if len(cliques_containing_e)==1 \
                and (cliques_containing_e[0] not in mandatory_cliques):
            mandatory_cliques.append(cliques_containing_e[0])

    # Check to see if mandatory_cliques contains a clique cover
    edges_in_set_of_cliques = set([])
    for clique in mandatory_cliques:
        edges_in_clique = [(clique[i], clique[j]) 
                           for i in xrange(len(clique)) 
                           for j in xrange(i+1,len(clique))]
        edges_in_set_of_cliques.update(set(edges_in_clique))
    if len(edges_in_set_of_cliques) == num_edges:
        if bound is None or len(mandatory_cliques) <= bound:
            return mandatory_cliques
        else:
            # There are too many cliques.  Return None to be
            # consistent with the documentation, even though we
            # actually know the clique cover number (and it is greater
            # than bound).
            return None

    max_cliques = [c for c in max_cliques if c not in mandatory_cliques]
    if bound==None:
        stopping_point=len(max_cliques)
    else:
        stopping_point=min(len(max_cliques),bound-len(mandatory_cliques))

    starting_point = max(1,ceil(num_edges / largest_clique_edges) \
                             - len(mandatory_cliques))
    for i in [starting_point..stopping_point]:
        for set_of_cliques in Combinations(max_cliques,i):
            edges_in_set_of_cliques = set([])
            for clique in set_of_cliques+mandatory_cliques:
                edges_in_clique = [(clique[i], clique[j]) 
                                   for i in xrange(len(clique)) 
                                   for j in xrange(i+1,len(clique))]
                edges_in_set_of_cliques.update(set(edges_in_clique))
            if len(edges_in_set_of_cliques) == num_edges:
                return set_of_cliques+mandatory_cliques
    return None


def edge_clique_cover_approximate(graph, bound=None):
    """
    Returns a decent (though not necessarily tight) upper bound for
    the clique cover number, which is the minimum number of cliques
    necessary to cover all of the edges in a graph

    INPUT:
        graph -- the graph to be examined
        bound -- this argument is ignored

    OUTPUT:
        a list of lists of vertices for each clique in the cover

    EXAMPLE:
        sage: sorted([sorted(e) for e in edge_clique_cover_approximate(graphs.PathGraph(4))])
        [[0, 1], [1, 2], [2, 3]]
    """
    free_edges=graph.edges(labels=False)
    vertices=graph.vertices() 
    vertices.sort(key=graph.degree) # sort the vertices by degree
    clique_vertices=[] #list of clique vertices
    while len(free_edges)>0:
        clique_edge=None
        while clique_edge is None: 
            # find an edge adjacent to a vertex of minimum degree
            v=vertices[0]
            edge_exists, edge=exists(free_edges, lambda edge: v in edge)
            if edge_exists:
                clique_edge=edge
            else:
                vertices.pop(0)

        edge_clique=find_edge_max_clique(clique_edge,graph)
        clique_vertices.append(edge_clique)
        edges_in_clique = graph.subgraph(edge_clique).edges(labels=False)
        free_edges= [e for e in free_edges if e not in edges_in_clique]
    return clique_vertices

def min_rank_by_bounds(graph, tests = ['precomputed', 'order', 'zero forcing', 'not path', 'forbidden minrank 2', 'not planar', 'not outer planar', 'clique cover', 'diameter']):
    """
    Return dictionaries giving the upper and lower bounds from running
    the specified tests.  If tests is not set, then all applicable
    tests are run.

    INPUT:
        graph -- the graph for which to find bounds

    OUTPUT: 
        a list of 2 dictionaries; the upper and lower bounds, respectively.

    EXAMPLE:  
        sage: g = Graph({0: [1,2,4,6,7], 1: [3,5,6,7,8], 2: [4,6,8], 3: [4,7,6], 4: [6], 5: [6,8,7]})
        sage: min_rank_by_bounds(g)
        ({'zero forcing': 4},
         {'clique cover': 9,
          'not outer planar': 6,
          'not path': 7,
          'not planar': 5,
          'order': 8})
        sage: min_rank_by_bounds(g, tests=['zero forcing', 'order', 'not path'])
        ({'zero forcing': 4}, {'not path': 7, 'order': 8})
    """
    if isinstance(tests, str):
        tests = [tests]

    order = graph.order()
    
    lower_bound = {}
    upper_bound = {}

    if 'precomputed' in tests:
        mr = get_mr_from_list(graph)
        if mr is not False:
            lower_bound['precomputed'] = mr
            upper_bound['precomputed'] = mr

    if 'order' in tests:
        upper_bound['order'] = order - 1

    if 'zero forcing' in tests:
        lower_bound['zero forcing'] = order - find_Z(graph)
        # Check if graph is a tree.  
        # If yes, then the ZFS will determine minimum rank.
        if graph.is_tree():
            upper_bound['zero forcing (tree)'] = lower_bound['zero forcing']

    if 'not path' in tests:
        if graph.diameter() < order - 1:
            upper_bound['not path'] = order - 2

    if 'forbidden minrank 2' in tests:
        if has_forbidden_induced_subgraph(graph):
	    lower_bound['forbidden minrank 2'] = 3
	else:
            upper_bound['forbidden minrank 2'] = 2

    if 'diameter' in tests:
	lower_bound['diameter'] = graph.diameter()

    if 'not planar' in tests:
        # Old versions of Sage assume that planar testing does not
        # have vertices of degree zero.  We can delete vertices of
        # degree zero without affecting the planarity.
        h = graph.copy()
        h.delete_vertices([v for v in h.vertices() if h.degree(v) == 0])
        if h.order()>0 and h.is_planar() is False:
            upper_bound['not planar'] = order - 4
    
    if 'not outer planar' in tests:
        if is_outerplanar(graph) is False:
            upper_bound['not outer planar'] = order - 3

    if 'clique cover' in tests:
        upper_bound['clique cover'] = len(edge_clique_cover(graph))
        
    return (lower_bound, upper_bound)


def find_cut_vertex(graph):
    """
    Return a "good" cut-vertex for a graph if it exists; otherwise,
    returns False.

    INPUT:
        graph -- the graph on which to find a cut-vertex

    OUTPUT:
        a cut-vertex (if one exists) that either results in components
        of order less than 7 or a minimum of the maximum component
        order; otherwise False

    EXAMPLES:
        sage: find_cut_vertex(graphs.PathGraph(3))
        1
        sage: find_cut_vertex(graphs.PathGraph(20))
        9
        sage: [find_cut_vertex(graphs.PathGraph(i)) for i in [1..20]]
        [False, False, 1, 1, 1, 1, 1, 1, 2, 3, 4, 5, 6, 6, 7, 7, 8, 8, 9, 9]
        sage: find_cut_vertex(graphs.CompleteGraph(3))
        False
    """

    vertices=graph.vertices()
    graph_cc_num=graph.connected_components_number()
    graph_order=graph.order()
    
    #this will hold the "best" cut-vertex and the order of the largest
    #connected component after deletion
    best_v=(False,graph_order)

    #checks each vertex and determines the best one
    for v in vertices:
        g=graph.copy()
        g.delete_vertex(v)
        g_cc = g.connected_components()
        if len(g_cc)>graph_cc_num:
            # We have a cut-vertex
            max_order = max(len(c) for c in g_cc)

            if max_order<7:
                return v
            if max_order<best_v[1]:
                best_v=(v,max_order)
            
    return best_v[0]

        
def find_rank_spread(vertex, graph):
    """
    Returns the exact rank spread for a graph and a vertex (i.e.,
    mr(G)-mr(G-v)) if the minimum ranks of both the graph and the
    graph without the vertex can be calculated using the minrank_bounds
    function.

    INPUT:
        vertex -- the vertex
        graph -- the graph

    OUTPUT:
        the rank spread and mr(graph-vertex) if both can be calculated
        using the minrank_bounds program, or False and False if either cannot
        be calculated exactly.
    
    EXAMPLES:
        sage: find_rank_spread(2, Graph({0:[1,2,3],1:[2,3],2:[3]}))
        (0, 1)
        sage: g = Graph({0:[1,2,4,6,7],1:[3,5,6,7,8], 2:[4,6,8],3:[4,7,6],4:[6],5:[6,8,7]})
        sage: find_rank_spread(2,g)
        (False, False)
    """
    subgraph=graph.copy()
    subgraph.delete_vertex(vertex)
    graph_bounds=minrank_bounds(graph)
    if graph_bounds[0]==graph_bounds[1]:
        # We have an actual min rank for graph
        subgraph_bounds=minrank_bounds(subgraph)
        if subgraph_bounds[0]==subgraph_bounds[1]:
            # We have an actual min rank for the subgraph
            return graph_bounds[0]-subgraph_bounds[0], subgraph_bounds[0]
    return False,False

                         
def cut_vertex_connected_graph_mr(c_vertex,graph):
    """
    Given a cut vertex and a graph, attempt to calculate the minimum
    rank of the graph by applying the cut vertex method to the graph
    and vertex.

    INPUT:    
        c_vertex -- the cut vertex
        graph -- the graph in which the cut vertex is contained

    OUTPUT:    
        a list of length 2 with the minimum rank as all entries, if
        the minimum rank can be calculated in this way

        False if the minimum rank cannot be calculated in this way

    EXAMPLE:
        sage: cut_vertex_connected_graph_mr(0,Graph({0:[1,2,3],2:[3]}))
        (2, 2)
        sage: cut_vertex_connected_graph_mr(2,Graph({0:[1,2,3],2:[3]}))
        Traceback (most recent call last):
        ...
        ValueError: Supplied vertex is not a cut vertex
    """
    g=graph.copy()
    if g.is_connected() is False:#this should never happen
        raise ValueError, "Graph is not connected"
    if c_vertex not in graph.vertices():#again, should never happen
        raise ValueError, "Supplied vertex is not in the graph"
    g.delete_vertex(c_vertex)
    subgraphs=g.connected_components_subgraphs()

    if len(subgraphs) <= 1: # c_vertex is not a cut-vertex
        raise ValueError, "Supplied vertex is not a cut vertex"
    
    index=0
    rank_spread=0
    subgraph_mr_sum=0
    for subgraph in subgraphs:
        subgraph_with_v = graph.subgraph(subgraph.vertices()+[c_vertex])
        new_rank_spread, subgraph_mr = \
            find_rank_spread(c_vertex, subgraph_with_v)
        if new_rank_spread is False:
            return False
        else:
            rank_spread += new_rank_spread
            subgraph_mr_sum += subgraph_mr

    rank_spread = min(rank_spread, 2)
    return subgraph_mr_sum+rank_spread, subgraph_mr_sum+rank_spread


def minrank_bounds(graph, all_bounds=False, tests=['precomputed', 'order', 'zero forcing', 'not path', 'forbidden minrank 2', 'not planar', 'not outer planar', 'clique cover', 'cut vertex', 'disconnected', 'diameter']):
    """
    Find lower and upper bounds for the minimum rank of a graph.  If
    all_bounds is False, then only return the best lower and upper
    bounds.  If True, return two dictionaries giving all applicable
    lower bounds and upper bounds, respectively.

    INPUT:
        graph -- the graph whose minimum rank is bounded
        
        all_bounds -- if False, then only return the best lower and upper
            bounds.  If True, return dictionaries giving all applicable lower
            bounds and upper bounds.

        tests -- a list of tests to get bounds.  Possible values are
            'precomputed', 'order', 'zero forcing', 'not path', 
            'no forbidden', 'not planar', 'not outer planar', 'clique cover',
            'cut vertex', 'disconnected'

    OUTPUT:
        the lower and upper bounds for the minimum rank, in that order


    EXAMPLES:
        sage: minrank_bounds(graphs.CompleteGraph(3))
        (1, 1)
        sage: minrank_bounds(graphs.CompleteGraph(3), all_bounds=True)
        ({'precomputed': 1, 'rank': 0, 'zero forcing': 1},
        {'clique cover': 1,
        'no forbidden': 2,
        'not path': 1,
        'order': 2,
        'precomputed': 1,
        'rank': 3})
        sage: minrank_bounds(graphs.PathGraph(4), all_bounds=True)
        ({'cut vertex (1)': 3, 'precomputed': 3, 'rank': 0, 'zero forcing': 3},
        {'clique cover': 3,
        'cut vertex (1)': 3,
        'order': 3,
        'precomputed': 3,
        'rank': 4,
        'zero forcing (tree)': 3})
        sage: minrank_bounds(graphs.PathGraph(4), all_bounds=True, tests=['order', 'zero forcing'])
        ({'rank': 0, 'zero forcing': 3},
        {'order': 3, 'rank': 4, 'zero forcing (tree)': 3})
        sage: minrank_bounds(graphs.HeawoodGraph())
        (8, 10)
        sage: minrank_bounds(graphs.HeawoodGraph(), all_bounds=True)
        ({'rank': 0, 'zero forcing': 8},
        {'clique cover': 21,
        'not outer planar': 11,
        'not path': 12,
        'not planar': 10,
        'order': 13,
        'rank': 14})
        sage: minrank_bounds(graphs.PetersenGraph())
        (5, 6)
        sage: minrank_bounds(graphs.PetersenGraph(), all_bounds=True)
        ({'rank': 0, 'zero forcing': 5},
        {'clique cover': 15,
        'not outer planar': 7,
        'not path': 8,
        'not planar': 6,
        'order': 9,
        'rank': 10})
    """
    if isinstance(tests, str):
        tests = [tests]

    possible_tests = set(['precomputed', 'order', 'zero forcing', 'not path', 'forbidden minrank 2', 'not planar', 'not outer planar', 'clique cover', 'cut vertex', 'disconnected', 'diameter'])
    # Check tests
    unknown_tests = set(tests).difference(possible_tests)
    if len(unknown_tests)>0:
        print "Unknown tests specified: ", list(unknown_tests)

    g=graph.copy()
    lower_bound = {'rank': 0}
    upper_bound = {'rank': g.order()}

    if g.is_connected():
        bounds = min_rank_by_bounds(graph, tests=tests)
        lower_bound.update(bounds[0])
        upper_bound.update(bounds[1])

        # Try finding a cut vertex
        if 'cut vertex' in tests:
            c_vertex=find_cut_vertex(g)
            if c_vertex is not False:
                cut_vertex_bounds = cut_vertex_connected_graph_mr(c_vertex,graph)
                if cut_vertex_bounds is not False:
                    lower_bound['cut vertex (%s)'%(c_vertex,)] = cut_vertex_bounds[0]
                    upper_bound['cut vertex (%s)'%(c_vertex,)] = cut_vertex_bounds[1]
    else:
        if 'disconnected' in tests:
            connected_components = g.connected_components_subgraphs()
            lower_bound['disconnected'] = 0
            upper_bound['disconnected'] = 0
            for component in g.connected_components_subgraphs():
                sub_bound = minrank_bounds(component, tests=tests)
                lower_bound['disconnected'] += sub_bound[0]
                upper_bound['disconnected'] += sub_bound[1]

    # Make sure that the lower bound is not greater than the upper bound
    if max(lower_bound.values()) > min(upper_bound.values()):
        raise StandardError, """
Best lower bound is greater than best upper bound; something is wrong:
lower bounds: %s
upper bounds: %s"""%(lower_bound,upper_bound)

    if all_bounds is True:
        return lower_bound, upper_bound
    else:
        # Return the best lower and upper bounds
        return max(lower_bound.values()), min(upper_bound.values())



def write_spreadsheet( graph_numbers, filename):
    """
    Write a spreadsheet of minimum ranks and bounds for Atlas of Graphs graphs.

    INPUT:
        graph_numbers -- a list of Atlas of Graph numbers
        
        filename -- The filename for the spreadsheet.

    OUTPUT:
        the file is written (or overwritten if it already exists)

    EXAMPLES: 
        To write a spreadsheet of graph minimum ranks for all graphs
        in the Atlas of Graphs similar to the one that comes with this
        program, first comment out the 'precomputed' test in the
        min_rank_by_bounds function, then do

        sage: write_spreadsheet([1..1252], "minrank.csv")
    """
    import csv
    f = file(filename, "wb")
    fieldnames = ('Atlas number','Order', 'Size','Minimum rank by program',
                  'Program lower bound','Program upper bound', 'Connected', 
                  'Zero forcing LB', 'Diameter LB', 'Clique cover UB',
                  'Not planar UB','Not outer planar UB','Not path UB',
                  'Induced subgraph (minrank 2)', 'Cut vertex','Tree')

    writer = csv.DictWriter(f, fieldnames=fieldnames)
    headers = {}
    for n in fieldnames:
        headers[n] = n
    writer.writerow(headers)

    for graph_number in graph_numbers:
        g = atlas_graphs[graph_number].copy()
        lower,upper = minrank_bounds(g, all_bounds=True)
        row = {}
        row['Atlas number'] = int(graph_number)
        row['Order'] = g.order()
        row['Size'] = g.size()

        if max(lower.values()) == min(upper.values()):
            row['Minimum rank by program'] = max(lower.values())
        else:
            row['Minimum rank by program'] = None

        row['Program lower bound'] = max(lower.values())
        row['Program upper bound'] = min(upper.values())

        row['Connected'] = g.is_connected()

        row['Zero forcing LB'] = lower.get('zero forcing', None)
        row['Diameter LB'] = lower.get('diameter', None)

        row['Clique cover UB'] = upper.get('clique cover', None)
        row['Not planar UB'] = upper.get('not planar', None)
        row['Not outer planar UB'] = upper.get('not outer planar', None)
        row['Not path UB'] = upper.get('not path', None)

        # If we have a forbidden minrank 2 lower bound, then we have a
        # forbidden subgraph
        if 'forbidden minrank 2' in lower:
            row['Induced subgraph (minrank 2)'] = True
        elif 'forbidden minrank 2' in upper:
            row['Induced subgraph (minrank 2)'] = False
        else:
            row['Induced subgraph (minrank 2)'] = None
        

        row['Cut vertex'] = False
        for key in lower.keys():
            if key.startswith('cut vertex'):
                row['Cut vertex'] = True
                break

        if 'zero forcing (tree)' in upper:
            row['Tree'] = True
        else:
            row['Tree'] = False

        writer.writerow(row)

    f.close() 
\end{lstlisting}

\end{document}